\documentclass[final]{siamltex}
\usepackage{amsmath,amssymb,latexsym}
\usepackage[final]{graphicx}
\newcommand{\field}[1]{\mathbb{#1}}
\newcommand{\R}{\field{R}}
\newcommand{\Z}{\field{Z}}
\newcommand{\N}{\field{N}}
\newcommand{\C}{\field{C}}

\newcommand{\F}{{\mathcal F}}

\newcommand{\tr}{\operatorname{\rm tr}}

\title{Computation of the entropy of polynomials orthogonal
on an interval}

\author{V.\ Buyarov\thanks{Moscow State University (RUSSIA)},
email: {\tt  BUYAROV@nw.math.msu.su}
 \and
J.\ S.\ Dehesa\thanks{Instituto Carlos I de F{\'\i}sica Te{\'o}rica y
Computacional, Granada University (SPAIN)}, email: {\tt
dehesa@ugr.es} \and A.\ Mart{\'\i}nez-Finkelshtein\thanks{University of
Almer{\'\i}a and Instituto Carlos I de F{\'\i}sica Te{\'o}rica y Computacional,
Granada University (SPAIN)}, corresponding author, email: {\tt
andrei@ual.es}  \and J.\ S{\'a}nchez-Lara\thanks{University of Almer{\'\i}a
(SPAIN)}, email: {\tt jlara@ual.es}
 }%
\date{\today}
\begin{document}

\maketitle

\begin{abstract}
We give an effective method to compute the entropy for polynomials
orthogonal on a segment of the real axis that uses as input data
only the coefficients of the recurrence relation satisfied by
these polynomials. This algorithm is based on a series expression
for the mutual energy of two probability measures naturally
connected with the polynomials. The particular case of Gegenbauer
polynomials is analyzed in detail. These results are applied also
to the computation of the entropy of spherical harmonics,
important for the study of the entropic uncertainty relations as
well as the spatial complexity of physical systems in central
potentials.
\end{abstract}



\pagestyle{myheadings} \thispagestyle{plain} \markboth{BUYAROV,
DEHESA, MART{\'I}NEZ-FINKELSHTEIN AND S{\'A}NCHEZ-LARA}{ENTROPY
COMPUTATION FOR ORTHOGONAL POLYNOMIALS}

\section{Introduction}

The concept of information entropy in its continuous or discrete
form has proved to be very fertile in numerous scientific branches
because of its flexibility and multiple meanings
\cite{Grad61,Grandy97,Ohya93,Zurek90}. Indeed, it is used as a
measure of disorder in thermodynamics \cite{Szilard29}, as a
measure of uncertainty in statistical mechanics \cite{Jaynes57} as
well as in classical and quantum information science
\cite{Holevo98,Nielsen00}, as a measure of diversity in ecological
structures, and as a criterion of classification of races and
species in population dynamics \cite{Karlin81}, among others.

In quantum mechanics, the uncertainty in the localization of a
particle in ordinary space is quantitatively measured by the
so-called position information entropy \cite{Bialynicki-Birula:75}
\begin{equation}\label{uno}
S_{\rho} = - \int \rho (\vec{r}) \ln \rho (\vec{r}) d \vec{r},
\end{equation}
in a better and more convenient way than the Heisenberg's standard
deviation of the quantum-mechanical probability density $\rho
(\vec{r}) = \left| \psi (\vec{r}) \right|^{2}$, where $\psi
(\vec{r})$ is the wavefunction of its dynamical state. Similarly,
the uncertainty in predicting the momentum of the particle is
measured by the momentum information entropy $S_{\gamma}$ of the
density $\gamma (\vec{p}) = \left| \tilde{\psi} (\vec{p})
\right|^{2}$, where the Fourier transform $\tilde{\psi} (\vec{p})$
of $\psi (\vec{r})$ is the wavefunction of the same state in the
dual, conjugate or momentum space. These two quantities describe
best \cite{Shannon48,Shannon49} the extent or spread of the
position and momentum probability densities, respectively.
Moreover, both entropies may decrease without bound when the
corresponding density becomes more concentrated, i.e. when
information in the associated space decreases. However, the
entropy sum is bounded from below
\cite{Beckner75,Bialynicki-Birula:75}
\[
S_{\rho} + S_{\gamma} \geq D (1 + \ln \pi),
\]
where $D$ is the dimensionality of the space (i.e. $D = 3$ for
ordinary space). It expresses the impossibility to have a complete
information of the position and momentum of the particle
simultaneously. This is the so-called entropic uncertainty
relation, which is a stronger version of the celebrated
Heisenberg's uncertainty principle, a fundamental law of nature.
This fact and the effective implementation of the density
functional theory of complex many-electron systems \cite{Parr89},
which uses the single-particle density as the basic variable, are
responsible for the fact that the study of the entropy has become
an ubiquitous tool in some areas (e.g.\ atomic and molecular
physics, condensed matter theories). For instance, several maximum
entropy methods based on the position and momentum entropies, as
well as on their sum, have been developed \cite{Gadre85,Jaynes57}
and widely used \cite{Angulo99,Gadre02,Zarzo97} for determination
of macroscopic quantities of natural systems. Nevertheless, the
lack of a general theoretical methodology and of accurate
numerical algorithms of computation these information entropies
still prevents this approach from being more widely used.

The exact determination of the information entropies of complex
many-particle systems is a formidable task. Only recently a small
progress has been achieved using the theory of special functions,
which in some cases allows to find closed formulas for the
information entropies of the simplest 1-dimensional
single-particle systems and the three-dimensional systems of
particles moving in a central or spherically symmetric potential.
For these systems the wavefunctions are controlled by some
classical orthogonal polynomials (such as Gegenbauer, Laguerre or
Hermite), and the determination of the corresponding information
entropies boils down naturally to the computation of entropic
functionals for sequences of orthogonal polynomials (cf.\
\cite{Angulo99, Nielsen00,Yanez:99}; a state-of-the art of this
topic up to 2001 is given in \cite{Dehesa:01}).

Namely, for a positive unit Borel measure  $\mu$ on $[-1,1]$,  let
\begin{equation}\label{polynomials}
p_n(x)=\gamma _n \prod_{j=1}^{n} \left ( x-\zeta_j^{(n)}   \right
)\,, \quad \gamma _n>0\,,
\end{equation}
be the corresponding \emph{orthonormal} polynomials:
\begin{equation*} \label{enint}
\int  p_n(x) p_m(x) \,   d \mu (x) = \delta_{mn}, \quad m, n \in
\Z_+\,.
\end{equation*}
We define the entropy of the polynomials $p_n(x)$ as
\begin{equation} \label{spg}
E_n=E_n(\mu ) = - \int p^2_n(x)  \, \ln \big( p^2_n(x) \big) \,  d
\mu (x) \,.
\end{equation}

The presence of these integrals raises two questions. One is the
study of their asymptotic behavior when $n \to \infty$, which has
a special interest in the analysis of the highly-excited (Rydberg)
states of numerous quantum-mechanical systems of hydrogenic-type
\cite{Yanez:94}. In this sense there have been important
contributions in the last few years
\cite{Aptekarev:95,Buyarov:99a,MR1790053,Dehesa:98,Lara02}; for a
detailed review, see \cite{Dehesa:01}. A totally different problem
is the  explicit computation of (\ref{spg}) for every fixed $n$
(up to a certain degree). Observe that a naive numerical
evaluation of these functionals by means of quadratures is not
convenient: since all the zeros of $p_n$ belong to the interval of
orthogonality, the increasing amount of integrable singularities
spoils any attempt to achieve reasonable accuracy even for rather
small $n$.

In this paper we present some theoretical results (Section
\ref{sec:series}), which allowed us to develop an algorithm for an
effective and accurate numerical computation of the entropy of
polynomials orthogonal on a segment of the real axis from the
coefficients of the three-term recurrence relation which they
satisfy (Section \ref{sec:computation}). In Section
\ref{section:Gegenb}, we study in detail the case of Gegenbauer
polynomials because of its own interest (as a very
``representative'' class of polynomials (\ref{polynomials})) and
because of their numerous applications; for instance, these
polynomials control the angular component of the wavefunctions of
single-particle systems in central potentials (cf.\
\cite{Dehesa:01}). In Section \ref{sec:experiments}, the results
of several numerical experiments are discussed, illustrating both
the accuracy and efficiency of the algorithm proposed here, and
comparing it with other computing strategies used so far. Finally,
the entropy of the spherical harmonics, which measures the spatial
complexity of single-particle systems and physical systems with
central potentials, is computed using the known relationship
between spherical harmonics and Gegenbauer polynomials (Section
\ref{sec:computationSH}).

\section{Series representation of the entropy} \label{sec:series}

The entropic functionals (\ref{spg}) can be restated in terms of
the logarithmic potential theory. If $\mu$ and $\nu$ are Borel
(generally speaking, real signed) measures on $\C$, we denote by
$$ V(z; \mu) = - \int \ln |z-t|\, d \mu (t)
$$ the logarithmic potential of $\mu$, and define the following functionals:
$$ I[\nu, \mu] = \int
V(z;\nu) \, d\mu(z)= - \iint \ln |z-t|\, d \nu (t) \, d\mu(z) \,,
$$ is the mutual energy of $\mu$ and $\nu$; $I(\mu)=I[\mu, \mu]$
is the logarithmic energy of $\mu $, and
$$ R[\nu,\mu] = -\int \ln \left(\frac{d \nu}{d \mu}\right) \, d
\nu
$$ is the relative entropy or the Kullback-Leibler information of
$\nu$ and $\mu$. From the Jensen inequality it immediately follows
that if both $\mu$ and $\nu$ are positive unit measures, then $
R[\nu,\mu] \leq 0$.

With the sequence of polynomials (\ref{polynomials}) we can
associate naturally two sequences of probability measures on
$[-1,1]$:
 \begin{equation}\label{zerocounting}
\lambda_n =\frac{1}{n}\,  \sum_{j=1}^n \delta_{\zeta_j^{(n)}}
\qquad \text{and} \qquad
d \nu_n(x)=p^2_n(x)\, d \mu(x)\,.
\end{equation}
Both measures are standard objects of study in the analytic theory
of orthogonal polynomials. For instance, the normalized zero
counting measure $\lambda _n$ is closely connected with the $n$-th
root asymptotics of $p_n$, and as was shown by Rakhmanov in his
pioneering work \cite{Rakhmanov:77}, $\nu_n$ is associated with
the behavior of the ratio $p_{n+1}/p_n$ as $n \to \infty$.

With the notations introduced above the entropy (\ref{spg}) is
equivalently rewritten as
\begin{equation}\label{kullback1}
E_n = R[\nu_n,\mu]\leq 0\,,
\end{equation}
or as
\begin{equation}\label{kullback2}
E_n =  - 2 \ln \gamma _n +  2\sum_{j=1}^n V(\zeta_j^{(n)};\nu_n)
=- 2 \ln \gamma _n + 2n \, I[\lambda _n, \nu _n]\,.
\end{equation}

A standard procedure for computation of $E_n$ employed so far
(once the quadratures begin to fail), is based on formula
(\ref{kullback2}). As it was shown in \cite{Dehesa:97}, the
potential $V(x;\nu_n)$ oscillates on $[-1,1] $ around the value
$\ln 2$, which is the Robin (or extremal) constant of this
interval. Moreover, the zeros $\zeta_j^{(n)}$ are points of local
minima for the potential $V(x;\nu_n)$, hence in order to compute
$E_n(\mu )$ we need to sum up the values of the logarithmic
potential $V(z; \nu_n)$ at its local minima. Nevertheless, this
procedure assumes an explicit computation of the zeros of $p_n$,
with the drawback of the well-known instability and computational
cost of this task as $n$ grows large.

An exception in this sense is an algorithm, proposed in
\cite{MR1790053} for numerical computation of the entropy of
Gegenbauer polynomials with integer parameters. It also involves
finding zeros of certain polynomials generated recursively, but
unlike in (\ref{kullback2}), the number of the zeros depends only
on the parameter of the polynomial and not on its degree.

In this paper we propose a totally different approach to the
computation of the entropic integrals (\ref{spg}), more in the
spirit of standard numerical algorithms for orthogonal
polynomials: it uses only the coefficients of the recurrence
relation satisfied by these polynomials as the input data. It does
not involve a solution of any nonlinear equation, and it can be
carried out by performing matrix multiplication of structured
(essentially, sparse) matrices. The algorithm is based on the
formula contained in Theorem \ref{theoremMainNew}, which seems to
be new.

We denote by $L^1_\mu$ the class of measurable functions on
$[-1,1]$, absolutely integrable with respect to $\mu$. Let $T_k
(x)=\cos (k \arccos x)$ denote as usual the Chebyshev polynomials
of the first kind. With the two measures introduced in
(\ref{zerocounting}) we define the double sequence of generalized
moments
\begin{equation}\label{defMandC}
    c_{k,n}=\int T_k(x)\, d \lambda _n(x)\,, \quad m_{k,n}=\int T_k(x)\, d \nu
    _n(x)\,, \quad k, n \geq 0\,.
\end{equation}
Obviously, $|c_{k,n}|\leq 1$ and $|m_{k,n}|\leq 1$ for all values
of $k$ and $n$. One of the main results of this paper is the
following
\begin{theorem} \label{theoremMainNew}
Assume that $\mu$ is a unit Borel measure on $[-1,1]$. For $n \in
\N$, let $p_n$ be the corresponding orthonormal polynomial, and
the unit measures $\lambda _n$ and $ \nu _n$ as defined in
(\ref{zerocounting}). Then for their mutual energy the following
formula holds:
\begin{equation}\label{series for energyNew}
    I[\lambda _n,\nu _n]=\ln 2+ 2 \, \sum_{k=1}^\infty
    \frac{c_{k,n}\,
    m_{k,n}}{k}\,,
\end{equation}
where the series in the right hand side is convergent.

Moreover, if we denote
\begin{equation}\label{condSuf1New}
M_n:=  \sup_{x \in [-1,1] } \, \int_{-1}^1 \left| \frac{p_n^2(x)-
p_n^2(t)}{x-t} \right|
 \, d\mu(x) <+\infty\,,
\end{equation}
then, for $N \in \N$ we have
\begin{equation}\label{remainderNew}
\left| I[\lambda _n,\nu _n]-\ln 2- 2 \, \sum_{k=1}^N \frac{c_{k,n}
m_{k,n}}{k}  \right|\leq \frac{4 M_n }{ N+1}\,.
\end{equation}
\end{theorem}
\begin{proof}
Following \cite{Dehesa:97}, we use the Fourier series of the
logarithm \cite[formula 1.514]{Gradshtein95},
\begin{equation}\label{logexpandCos}
-\ln |1-e^{i\varphi }|= \sum_{k=1}^\infty \frac{\cos k \varphi
}{k}\,,
\end{equation}
valid for almost all $\varphi \in [0,\pi]$ (see e.g.\
\cite[Theorem 15.2]{Champeney87}), which yields a representation
of the logarithmic kernel
\begin{equation}\label{logexpand}
-\ln |x-t|=\ln 2 + 2 \sum_{k=1}^\infty \frac{1}{k}\, T_k(x)
T_k(t)\,,
\end{equation}
where for every $x \in [-1,1]$ the series (in $t$) converges
almost everywhere in $[-1,1]$.

On the other hand, by the recurrence relation
\begin{equation}\label{recurrenceChebA}
T_{k+1}(x)=2 x\,  T_{k}(x)-T_{k-1}(x)\,, \quad T_{0}(x)=1, \quad
T_{1}(x)=x\,,
\end{equation}
we have that
$$
2(t-x)  T_k(x) T_k(t)= q_{k}(x,t)-q_{k-1}(x,t)\,, \quad
q_{k}(x,t)= T_{k+1}(t) T_k(x)-T_{k+1}(x) T_k(t)\,,
$$
from where for $N \in \N$,
$$
2(t-x)  \sum_{k=1}^N \frac{T_k(x) T_k(t)}{k}=
-(t-x)+\sum_{k=1}^{N-1} \frac{q_{k}(x,t)}{k (k+1)}+ \frac{
q_{N}(x,t)}{N}\,.
$$
Hence,
$$
\left| 2(t-x)  \sum_{k=1}^N \frac{T_k(x) T_k(t)}{k} \right| \leq
4\,, \quad t, x \in [-1,1]\,.
$$

Thus, if $f \in L_\mu^1$, we can apply Lebesgue dominated
convergence theorem to (\ref{logexpand}) in order to assert that
\begin{equation}\label{mainFormula1}
-\int_{-1}^1 f(t) (t-x)  \ln |t-x|\, d \mu(t) = \widehat f_0\, \ln
2  + 2 \sum_{k=1}^\infty \frac{\widehat f_k}{k}\, T_k(x) \,, \quad
\widehat f_k=\int_{-1}^1 f(t) (t-x) T_k(t)\, d \mu(t)\,,
\end{equation}
and the series in the right hand side is convergent. Furthermore,
we can estimate the remainder using that
$$
2(t-x)  \sum_{k=N+1}^\infty \frac{T_k(x) T_k(t)}{k}=
-\frac{q_{N}(x,t)}{N+1}+\sum_{k=N+1}^\infty \frac{q_{k}(x,t)}{k
(k+1)}\,,
$$
from where
$$
\left| 2(t-x)  \sum_{k=N+1}^\infty \frac{T_k(x) T_k(t)}{k} \right|
\leq \frac{4}{N+1}\,, \quad t, x \in [-1,1]\,,
$$
so that
\begin{equation}\label{mainFormulaX}
\left| 2 \sum_{k=N+1}^\infty \frac{\widehat f_k}{k}\, T_k(x)
\right| \leq \frac{4}{N+1}\, \int_{-1}^1 |f(t)| \, d \mu(t)\,.
\end{equation}

In particular, for $x \in [-1,1]$ we may take
$$
f(t,\cdot)=\frac{p_n^2(t)-p_n^2(x)}{t-x}\,,
$$
and formula (\ref{mainFormula1}) yields
\begin{equation*}\label{mainFormula2}
\begin{split}
-\int_{-1}^1 \big(p_n^2(t)-p_n^2(x)\big) \,  \ln |t-x|\, d \mu(t)
&= \\ \ln 2 \, \int_{-1}^1 \big(p_n^2(t)-p_n^2(x)\big) \, d \mu(t)
& + 2 \sum_{k=1}^\infty \frac{T_k(x)}{k}\,  \int_{-1}^1
\big(p_n^2(t)-p_n^2(x)\big) \, T_k(t)\, d \mu(t)\,,
\end{split}
\end{equation*}
which using the definitions introduced above can be rewritten as
\begin{equation}\label{potentials}
V(x;\nu_n)-p_n^2(x)V(x;\mu)=(1-p_n^2(x))\ln 2+ 2 \sum_{k=1}^\infty
\frac{T_k(x)}{k}\, \left(m_{k,n}-p_n^2(x)\, \int_{-1}^1  T_k(t)\,
d \mu(t) \right)\,.
\end{equation}
Evaluating (\ref{potentials}) at the zeros of $p_n$ we obtain
$$
V(\zeta_j^{(n)};\nu_n)= \ln 2+ 2 \sum_{k=1}^\infty
\frac{T_k(\zeta_j^{(n)})}{k}\,
 m_{k,n}\,, \quad j=1, 2, \dots, n\,.
$$
Summing up these expressions for $j=1, 2, \dots, n$, we arrive at
(\ref{series for energyNew}).

On the other hand, by (\ref{condSuf1New}) and
(\ref{mainFormulaX}),
$$
\left|2 \sum_{k=N+1}^\infty \frac{T_k(x)}{k}\,  \int_{-1}^1
\big(p_n^2(t)-p_n^2(x)\big) \, T_k(t)\, d \mu(t)\right| \leq
\frac{4}{N+1} \, \int_{-1}^1 \left | \frac{p_n^2(t)-p_n^2(x)}{t-x}
\right|\,d\mu(t)\leq \frac{4 M_n}{N+1} \,,
$$
from where
$$
\left| V(\zeta_j^{(n)};\nu_n)- \ln 2- 2 \sum_{k=1}^N
\frac{T_k(\zeta_j^{(n)})}{k}\,
 m_{k,n}\right| = \left|2 \sum_{k=N+1}^\infty
\frac{T_k(\zeta_j^{(n)})}{k}\,
 m_{k,n}  \right| \leq \frac{4 M_n}{N+1} \,, \quad j=1, 2, \dots,
 n\,,
$$
and the estimate  (\ref{remainderNew}) follows.
\end{proof}

\medskip

\begin{corollary}
With assumptions of Theorem \ref{theoremMainNew},
\begin{equation}\label{series for entropy}
E_n = - 2 \ln \left(    \frac{ \gamma _n}{  2^n }\right  ) + 4n \,
\sum_{k=1}^\infty \frac{c_{k,n} m_{k,n}}{k}\,,
\end{equation}
and the error after truncating the series at the $N$-th term is
bounded by the right hand side in (\ref{remainderNew}).
\end{corollary}

\medskip

\textsc{Remark:} As our example of Gegenbauer polynomials in
Section \ref{section:Gegenb} shows, the bound (\ref{remainderNew})
is usually too pessimistic.

\section{Effective computation of $E_n$}
\label{sec:computation}

Assume that we have as input data the coefficients of the
three-term recurrence relation, satisfied by the orthonormal
polynomials $p_n(x)=\gamma _n x^n+\dots$,
 \begin{equation}\label{3term}
x p_{n}(x)=a_{n+1} p_{n+1}(x)+ b_n p_n(x) + a_n p_{n-1}(x)\,,
\end{equation}
with $p_{-1}=0$ and $p_0(x)=1$. We form the infinite Jacobi matrix
$$
J=\begin{pmatrix} b_0 & a_1 & 0 & 0 &\dots  \\
a_1 & b_1  & a_2 & 0 & \dots  \\
0 & a_2 & b_2 & a_3 & \dots  \\
\vdots & \vdots & \ddots & \ddots & \ddots
\end{pmatrix}
$$
and let $J_{n}=J(1:n,1:n)$ denote its principal minor $n \times
n$. Here and in what follows we occasionally use the MATLAB type
notation to refer to elements of a matrix. Furthermore, $e_n$
stands for the infinite (column) vector whose $n$-th element is 1
and the rest is 0, and $\langle a, b \rangle=a^H b$ is the
standard scalar product in $l^2$ or $\R^k$ (the space where we
multiply vectors is always clear from the context).

The following are very well known facts:
\begin{proposition}
Let $p_n$ be the orthonormal polynomials (\ref{polynomials})
satisfying the recurrence relation (\ref{3term}). Then, with the
notation above, for $n \geq 1$,
\begin{enumerate}
 \item[i)] the zeros $\zeta_{j}^{(n)}$, $j=1, \dots, n$, of $p_n$ are eigenvalues of
$J_n$, and  $(p_0(\zeta_{j}^{(n)}), p_1(\zeta_{j}^{(n)}), \dots,
p_{n-1}(\zeta_{j}^{(n)})^T$ are corresponding eigenvectors. In
particular, for $m \in \N_0=\N \cup \{0\}$,
$$
\sum_{j=1}^n \left[\zeta_{j}^{(n)} \right]^m = \text{trace of }
[J_n]^m.
$$
 \item[ii)] If $f$ is a polynomial then
\begin{equation}\label{JbyJn}
\left \langle e_{n+1}, f(J) e_{n+1} \right \rangle =\int_{-1}^1
f(x) p_n^2(x) w(x)\, dx\,.
\end{equation}

 \item[iii)] For $r\geq n+1$ and $m \in \N_0$,
\begin{equation}\label{JbyJn2}
\left \langle e_{n+1}, [J_r]^m e_{n+1} \right \rangle
=\sum_{j=1}^r \Lambda_j^{(r)} \, \left[\zeta_{j}^{(r)}\right]^m\,
p_n^2\left(\zeta_{j}^{(r)}\right) \,,
\end{equation}
where $\Lambda_j^{(r)}$ are the Cotes-Christoffel numbers (Gauss
quadrature weights) given by
\begin{equation}\label{christoffel}
\Lambda_j^{(r)}=\left[\sum_{i=0}^{r-1}
p_i^2\left(\zeta_{j}^{(r)}\right) \right]^{-1}.
\end{equation}
 \item[iv)] The leading coefficient of $p_n$ satisfies $\gamma _n=(a_1 a_2 \dots a_n)^{-1}$.
\end{enumerate}
\end{proposition}
All these facts are classical (see e.g.\ \cite{szego:1975}).
Formula (\ref{JbyJn}) can be found for instance in \cite[\S
4.1.2]{vanassche:87}, where it is proved for $f \in C[-1,1]$.
Identity \eqref{JbyJn2} is a straightforward consequence of
\eqref{christoffel}, $i)$, and the spectral decomposition of the
finite selfadjoint matrix $J_r$.

For practical computation of the left hand side in \eqref{JbyJn}
we need the following result:
\begin{corollary} \label{cor:truncate}
If $n \in \N_0$, $m, r \in \N$ and $r \geq n +(m+1)/2$, then the
$(n+1,n+1)$ elements of $J^m$ and $(J_{r})^m$ coincide.
\end{corollary}
\begin{proof}
It is well known that Gauss quadrature formula
$$
\int_{-1}^1 f(x)   w(x)\, dx=\sum_{j=1}^r \Lambda_j^{(r)} \,
f\left(\zeta_{j}^{(r)}\right)
$$
is exact  if the degree of the polynomial $f$ is $\leq 2 r-1$. In
particular,
$$
\int_{-1}^1 x^m p_n^2 (x)   w(x)\, dx=\sum_{j=1}^r \Lambda_j^{(r)}
\, \left[\zeta_{j}^{(r)}\right]^m
p_n^2\left(\zeta_{j}^{(r)}\right)
$$
is exact if $m+ 2n \leq 2r-1$, and the statement follows from
\eqref{JbyJn}--\eqref{JbyJn2}.\footnote{This short and elegant
proof was suggested by one of the anonymous referees whose
contribution we gratefully acknowledge.}
\end{proof}

\medskip

Assume that computing the entropy by means of formula (\ref{series
for entropy}) we decide to truncate the series therein  at $k=N
\in \N$. Then as a consequence of the previous corollary, in the
right hand side of (\ref{JbyJn}) we may use $T_k(J_{r})$ instead
of $T_k(J)$, with $r = n+ 1+ [N/2]$. Moreover, in the right hand
side of (\ref{JbyJn}) it is sufficient to know only the $(n+1)$-th
column of $f(J_r)$. Recalling that Chebyshev polynomials satisfy
the recurrence relation (\ref{recurrenceChebA}), we can propose
the following algorithm, where $I_n$ stands for the $n \times n$
identity matrix:

\medskip

\begin{center}
\framebox{
\begin{minipage}{10cm}
\begin{center}
\textsc{Algorithm 1}
\end{center} \smallskip
\texttt{\begin{tabbing}
 \scriptsize{(i)} \phantom{1111} \= Compute \= $ \ln \left(\gamma _n/  2^n \right  )
 =-\ln \left [  \prod_{j=1}^n  \left( 2 a_j \right  )\right ]  $
 \=
recursively; \\
  \>choose $N \in \N $ at which truncate the series in
(\ref{series for entropy});\\
  \> take
$T_{0}(J_n)=I_n$, $T_{1}(J_n)=J_n$, and iterate \\
  \>  \> $T_{k}(J_n)=2  J_n\, T_{k-1}(J_n)-T_{k-2}(J_n)$, $k=2,
\dots, N$,
\\  \> computing 
\\
  \scriptsize{(ii)} \> \> $ c_{k,n}=  \tr\, T_{k}(J_n)/n$,   $ k=1, \dots, N$;
\\
  \> set $r = n +1+ [N/2]$; starting with $v_0=I_r(:,n+1)$ and
$v_1=J_r(:,n+1)$, \\
\> iterate by the recurrence \\
\> \> $v_{k}=2 J_r\, v_{k-1}-v_{k-2}$, $k=2, \dots, N$, \\  \>
computing 
\\
 \scriptsize{(iii)}  \> \> $m_{k,n}=v_{k}(n+1)$, $ k=1, \dots, N$; \\
  \> substitute the results of (i)-(iii)  in (\ref{series for entropy}),
terminating the series at $k=N$.
\end{tabbing}}
\end{minipage}
}
\end{center}

\medskip

Observe that this algorithm starts from the spectral data as the
only input, and performs multiplication of structured matrices,
without solving any kind of equation. This can be efficiently
implemented, for instance, using the known algorithms for sparse
matrix multiplication.

On the other hand, in order to satisfy conditions of Corollary
\ref{cor:truncate} it is necessary to know in advance the
truncation term $N$ for the series in (\ref{series for entropy}),
for which we need an a priori bound on the error. The bound in
(\ref{remainderNew}) can be used, but usually it overestimates the
error yielding values of $N$ much larger than needed. In Section
\ref{section:Gegenb}, we discuss the selection of the truncation
term  in the particular case of Gegenbauer polynomials.

%

\section{Entropy computation for Gegenbauer polynomials} \label{section:Gegenb}

In this Section we test our approach on the computation of the
entropy of the Gegenbauer polynomials.

For $\lambda >-1/2$ let
\begin{equation}\label{normalizationC}
c_\lambda =\frac{\Gamma(\lambda +1)}{\sqrt{\pi}\, \Gamma(\lambda
+1/2)}\,.
\end{equation}
It is easy to verify that $ w^\lambda (x)=c_\lambda
(1-x^2)^{\lambda -1/2} $ is a positive unit weight on $[-1,1]$.
Let $C_k^\lambda $ denote the Gegenbauer polynomial of degree $k$
and parameter $\lambda $, orthogonal  with respect to $w^\lambda$
on this interval, and normalized by the value at $x=1$,
\begin{equation}\label{standard}
C_k^\lambda (1)=\binom{k+2\lambda -1}{k}\,;
\end{equation}
this is a standard normalization, adopted for instance in
\cite{abramowitz/stegun:1972} and in \cite{szego:1975}.
Straightforward computation shows that
\begin{equation}\label{normalizationGegen}
G_k^\lambda(x)=\left( \frac{k!\, (k+\lambda) \Gamma(2
\lambda)}{\lambda \Gamma(k+2\lambda )} \right)^{1/2}\,
C_k^\lambda(x)=\gamma  _k^\lambda\,  x^k+\text{lower degree terms}
\end{equation}
are the Gegenbauer polynomials orthonormal with respect to
$w^\lambda (x)$.

In this Section we will use the superscript $\lambda $ in the
previously introduced notation when we want to make the dependence
on the parameter $\lambda $ explicit; for instance,
\begin{equation} \label{entropyGegenbDef}
E_n^{\lambda }= - \int_{-1}^1 \left (G_n^\lambda(x)\right)^2 \,
\ln \left (G_n^\lambda(x)\right)^2 w^\lambda (x) \,  d x \,.
\end{equation}

For the time being, only few explicit formulas for the entropy of
orthonormal Gegenbauer polynomials are known. Namely, for $\lambda
=0$ and $\lambda =1$ (Chebyshev polynomials of the first and
second kind, respectively) it is not difficult to prove (cf.\
\cite{Dehesa:97,Yanez:94}) that
\begin{equation} \label{ent1bis}
E_n^0=\log2-1,\qquad E_n^1=-\frac{n}{n+1}\,.
\end{equation}
Furthermore, case $\lambda=2$ was studied in \cite{Buyarov:97} and
\cite{MR1790053}, establishing that
\begin{equation}\label{ent2}
E_n^2=\log\left(\frac{n+3}{3(n+1)}\right)-\frac{n^3-5n^2-29n-27}{(n+1)(n+2)(n+3)}-
\frac{1}{n+2}\left(\frac{n+3}{n+1}\right)^{n+2}.
\end{equation}
This formula does not allow to expect a short and elegant
expression for all $E_n^\lambda $, even for integer values of
$\lambda $.

For $\lambda \in \N$, $\lambda \geq 2$, an alternative algorithm
has been proposed in \cite{MR1790053}; it expresses the entropy in
terms of the zeros of certain polynomials generated recursively,
whose number, $2\lambda -2$, depends only on the parameter of the
polynomial and not on its degree.
%

This approach (valid only for integer values of $\lambda $), is
more efficient than the direct computation of $E_n^{\lambda}$ by
quadrature when $n$ grows large, and allows also to find
constructively the first terms of the asymptotic expansion of the
entropy when $n \to \infty$ and $\lambda $ is fixed (cf.\
\cite{Aptekarev:95,MR1790053,Lara02}):
\begin{equation}\label{asympt_Gegen_entropy}
    E_n^\lambda =\mathcal{E}^\lambda_0+\frac{\mathcal{E}^\lambda_1}{n}+\dots\,,
    \quad \mathcal{E}^\lambda_0=-1-\ln \frac{\Gamma(2 \lambda )}{\Gamma(\lambda ) \Gamma(\lambda
    +1)}\,.
\end{equation}

In this Section we will apply the general approach, described in
Section \ref{sec:computation}, to the efficient computation of the
entropy $E_n^{\lambda}$. It is well known that the polynomials
$G_n^{\lambda}$ satisfy the three-term recurrence relation
$$
x \, G_{n}^\lambda (x)=a_{n+1} G^\lambda_{n+1}(x)+ a_n
G^\lambda_{n-1}(x)\,,
$$
where
$$
a_n=\frac{1}{2} \, \left[ \frac{n (n +2 \lambda -1) }{(n+\lambda
-1)(n+\lambda )} \right]^{1/2}>0\,.
$$
In particular, since the leading coefficient $\gamma
_n^\lambda=(a_1 a_2 \dots a_n)^{-1}$, we have
 \begin{equation}\label{gammaForGegen}
-2 \ln \left(\gamma  _n^\lambda/  2^n \right  )
 =\ln  \left[ \prod_{j=1}^n  \frac{j (j +2 \lambda -1) }{(j+\lambda
-1)(j+\lambda )} \right]\,.
\end{equation}
Consider the values defined in (\ref{defMandC}). By symmetry,
\begin{equation}\label{zeroForOdd}
c_{k,n}^{\lambda }=m_{k,n}^{\lambda }=0 \text{ for $k$ odd.}
\end{equation}
Observe that for $\lambda >-1/2$ the orthogonality weight
satisfies the assumptions of Theorem \ref{theoremMainNew}. Thus,
taking into account (\ref{gammaForGegen}), formula (\ref{series
for entropy}) for the Gegenbauer polynomials reads as
\begin{equation}\label{series_for_entropu_Gegen}
E_n^\lambda  = \ln  \left(\prod_{j=1}^n  \frac{j (j +2 \lambda -1)
}{(j+\lambda -1)(j+\lambda )} \right) + 2n \, \sum_{k=1}^{\infty}
\frac{c_{2k,n}^{\lambda } m_{2k,n}^{\lambda }}{k}\,, \quad \lambda
>-1/2\,, \quad n \in \N\,.
\end{equation}
\begin{proposition} \label{prop:terminating}
For $\lambda \in \N_0$, the series in
(\ref{series_for_entropu_Gegen}) is terminating at $k=n+\lambda $.
\end{proposition}
\begin{proof}
This is a straightforward consequence of the fact that in this
case $ w^\lambda (x)  (1-x^2)^{1/2}$ is a polynomial of degree
$\lambda $, and for $k> 2n+2\lambda $ we can use orthogonality of
$T_k$ in the definition of $m_{k,n}$ in (\ref{defMandC}).
\end{proof}

Taking advantage of (\ref{zeroForOdd}), we can use the recurrence
formula for even Chebyshev polynomials: $T_0(x)=1$,
$T_2(x)=2x^2-1$, and
\begin{equation}\label{recurrenceEvenCheb}
    T_{2k}(x)=2\, T_2(x)\, T_{2k-2}(x)-T_{2k-4}(x)\,, \quad k \geq
    2\,.
\end{equation}

Then Algorithm 1 described in Section \ref{sec:computation} takes
the following form  for the Gegenbauer polynomials:

\medskip

\begin{center}
\framebox{
\begin{minipage}{10cm}
\begin{center}
\textsc{Algorithm 2}
\end{center} \smallskip
\texttt{\begin{tabbing}
 \scriptsize{(i)} \phantom{1111} \= Find \= $ -2 \ln \left(\gamma _n^\lambda /  2^n \right  )
 $  \=  by (\ref{gammaForGegen}) recursively. \\
  \> Choose a value $N \geq n+\lambda $ where to truncate the
series in (\ref{series_for_entropu_Gegen}),  \\ \> \> unless
$\lambda \in
\N_0$; in such a case, $N=n+\lambda$.\\
  \> Take
$T_0(J_n)=I_n$, $\widehat J_n=T_2(J_n)=2J_n^2-I_n$, and iterate \\
  \>  \>  \phantom{111}  $T_{2k}(J_n)=2 \widehat J_n\, T_{2k-2}(J_n)-T_{2k-4}(J_n)$, $k=2,
\dots, N$,
\\  \>  \> computing 
\\
  \scriptsize{(ii)} \> \>  \phantom{111}  $c_{2k,n}^{\lambda }=  \tr\, T_{2k}(J_n)/ n$, $ k=1,
\dots, N$.
\\
  \> Set $r = n+ 1+ [N/2]$ and $\widehat J_r=T_2(J_r)=2J_r^2-I_r$. \\
  \> Starting with $v_0=I_r(:,n+1)$ and
$v_2=\widehat J_r(:,n+1)$, \\
\> \> iterate by the recurrence \\
\> \>  \phantom{111}  $v_{2k}=2 \widehat J_r\, v_{2k-2}-v_{2k-4}$,
$ k=2, \dots, N$, \\  \>  \>
computing 
\\
 \scriptsize{(iii)}  \> \>  \phantom{111}  $m_{2k,n}^{\lambda }=v_{2k}(n+1)$, $ k=1, \dots, N$. \\
  \> Substitute the results of (i)-(iii)  in
  (\ref{series_for_entropu_Gegen}),\\ \> \>
terminating the series at $k=N$.
\end{tabbing}}
\end{minipage}
}
\end{center}

\medskip

In the implementation it is convenient to use subroutines for
sparse matrix multiplication; as it was mentioned, the coefficient
in (\ref{gammaForGegen}) is better to compute recursively,
avoiding possible overflows.

In order to obtain an a priori bound for the error (and thus to
find the truncation term $N$) we can use the explicit expression
for the coefficients $m_{2k,n}^{\lambda }$. In \cite{Yanez:94},
$m_{2k,n}^{\lambda } $ were expressed in terms of Wilson
polynomials of degree $n$ and parameters depending on $\lambda $
and $k$. Nevertheless, the expression which appears there has
indeterminacies for integer values of $\lambda $, which is
inconvenient for evaluation. We are interested in an alternative
formula for $m_{2k,n}^{\lambda }$.

When a linearization formula of the type
$$
T_k(x) \, p_n(x) =\sum_{j=0}^{n+k} \ell_{j,k,n} p_j(x)\,,
$$
is available the coefficients $m_{k,n}$ in (\ref{defMandC}) can be
found observing that $ m_{k,n}=\ell_{n,k,n} $.  Nevertheless, we
were unable to find the explicit expression in the literature, and
we establish a formula for $m_{2k,n}^{\lambda }$ based on the
hypergeometric representation for the Chebyshev and Gegenbauer
polynomials.
\begin{theorem}
For the orthonormal Gegenbauer polynomials $G_{n}^\lambda$ the
coefficients $m_{2k,n}^{\lambda}$ defined in (\ref{defMandC})
satisfy
\begin{equation}\label{explicitm1}
m_{2k,n}^{\lambda}=\frac{n+\lambda}{n!}\sum_{j=0}^n(-1)^j
\binom{n}{j} \frac{(2\lambda+j)_n}{ j+\lambda}
\frac{(-j-\lambda)_k}{(j+\lambda+1)_k}\,, \quad k \geq 1\,.
\end{equation}
Alternatively, for $k>n+\lambda $,
\begin{equation}\label{mkn}
 m_{2k,n}^{\lambda}=-\sin(\pi \lambda )\, \frac{n+\lambda }{\pi \, n!}\,
\sum_{j=0}^n \binom{n}{j}\, \frac{(2\lambda +j)_n}{j+\lambda } \,
\Gamma^2 (j+\lambda +1) \, \frac{\Gamma(k-j-\lambda
)}{\Gamma(k+j+\lambda +1)}\,.
\end{equation}
In particular, for $\lambda \in \N_0$, $ m_{2k,n}=0$ for all $k
\geq n+\lambda +1$. Moreover, for any $\lambda >-1/2$ and
$k>n+\lambda$, we have $|m_{2k+2,n}| \leq |m_{2k,n}|$, and
\begin{equation}\label{asymptoticsMkn}
|m_{2k,n}| =O\left (  \frac{ 1  }{  k^{2\lambda +1}} \right )\,,
\quad   k\to +\infty\,.
\end{equation}
\end{theorem}
Here $(a)_k=\Gamma(a+k)/\Gamma(a)$ denotes the Pochhammer's
symbol.
\begin{proof}
The following hypergeometric representation for the Chebyshev and
Gegenbauer polynomials is well known (see e.g.\ \cite[formulas
8.942 and 8.932]{Jeffrey95}):
$$
T_{2k}(x)=~_2F_1\left(\left.\begin{array}{c}-2k,~2k\\1/2
\end{array}\right|\frac{1-x}{2}\right),\label{hyperRe1}\quad
C_n^{\lambda}(x)=\frac{(2\lambda)_n}{n!}~_2F_1\left(\left.\begin{array}{c}-n,~n+2\lambda\\
\lambda+1/2\end{array}\right|\frac{1-x}{2}\right).\label{hyperRe2}
$$
The quadratic transformation
\begin{equation*}
_2F_1\left(\left.\begin{array}{c}-2n,~2n+2\alpha+1\\
\alpha+1\end{array}\right|\frac{1-x}{2}\right)=
~_2F_1\left(\left.\begin{array}{c}-n,~n+\alpha+1/2\\
\alpha+1\end{array}\right|1-x^2\right)
\end{equation*}
yields
$$
T_{2k}(x)=~_2F_1\left(\left.\begin{array}{c}-k,~k\\1/2\end{array}\right|1-x^2\right),
\quad C_n^{\lambda}(x)=\frac{(2\lambda)_n}{n!}~ _2F_1\left(\left.
\begin{array}{c} -n/2,~n/2+\lambda \\ \lambda+1/2
\end{array}
\right|1-x^2\right),
$$
and Clausen's identity (cf.\ \cite[p.\ 116, problem
13]{Andrews:99}  or \cite[Ch.\ IV, section 4.3]{Erdelyi53a})
\begin{equation*}
\left(_2F_1\left(\left.\begin{array}{c}a,~b\\a+b+1/2\end{array}\right|x\right)\right)^2=~
_3F_2\left(\left.\begin{array}{c}2a,~2b,~a+b\\2a+2b,~a+b+1/2\end{array}\right|x\right)
\end{equation*}
gives the following representation:
\begin{equation*}
[C_n^{\lambda}(x)]^2=\left(\frac{(2\lambda)_n}{n!}\right)^2
~_3F_2\left(\left.\begin{array}{c}-n,~n+2\lambda,~\lambda\\2\lambda,~\lambda+1/2\end{array}
\right|1-x^2\right).
\end{equation*}
We use these formulas in order to compute the integral
\begin{align*}
D_{2k,n}&=\int_{-1}^1T_{2k}(x)~C_n^{\lambda}(x)^2~(1-x^2)^{\lambda-1/2}dx
=\left(\frac{(2\lambda)_n}{n!}\right)^2 \\
& \times
\sum_{i=0}^k\frac{(-k)_i~(k)_i}{(1/2)_i~i!}\left(\sum_{j=0}^n
\frac{(-n)_j~(n+2\lambda)_j~(\lambda)_j}{(2\lambda)_j
(\lambda+1/2)_j~j!}~\int_{-1}^1(1-x^2)^{i+j+\lambda-1/2}dx\right)\\
&=\sqrt{\pi}\left(\frac{(2\lambda)_n}{n!}\right)^2\sum_{i=0}^k\sum_{j=0}^n
\frac{(-k)_i~(k)_i}{(1/2)_i~i!}~
\frac{(-n)_j~(n+2\lambda)_j~(\lambda)_j}{(2\lambda)_j~(\lambda+1/2)_j~j!}
~\frac{\Gamma(i+j+\lambda+1/2)}{\Gamma(i+j+\lambda+1)}.
\end{align*}
Interchanging the order of summation we get
$$
D_{2k,n} =\sqrt{\pi}\left(\frac{(2\lambda)_n}{n!}\right)^2
\sum_{j=0}^n
~_3F_2\left(\left.\begin{array}{c}-k,~k,~j+\lambda+1/2\\1/2,~j+\lambda+1\end{array}
\right|1\right)
\frac{(-n)_j~(n+2\lambda)_j~(\lambda)_j}{(2\lambda)_j~(\lambda+1/2)_j~j!}
~\frac{\Gamma(j+\lambda+1/2)}{\Gamma(j+\lambda+1)}.
$$
Using the Pfaff-Saalschutz identity,
\begin{equation*}
~_3F_2\left(\left.\begin{array}{c}-k,~a,~b\\c,1+a+b-c-k\end{array}
\right|1\right)=\frac{(c-a)_k~(c-b)_k}{(c)_k~(c-a-b)_k},
\end{equation*}
the integral $D_{2k,n}$ becomes
\begin{equation*}
D_{2k,n}=\sqrt{\pi}\left(\frac{(2\lambda)_n}{n!}\right)^2\sum_{j=0}^n
\frac{(1/2-k)_k~(-j-\lambda)_k}{(1/2)_k~(-k-j-\lambda)_k}~
\frac{(-n)_j~(n+2\lambda)_j~(\lambda)_j}{(2\lambda)_j~(\lambda+1/2)_j~j!}
~\frac{\Gamma(j+\lambda+1/2)}{\Gamma(j+\lambda+1)}.
\end{equation*}
Taking into account the normalization factors in
(\ref{normalizationC}) and (\ref{normalizationGegen}) we obtain
\begin{align*}
m_{2k,n}^{\lambda}&=\frac{n!~(n+\lambda)}{\lambda~(2\lambda)_n}
~\frac{\Gamma(\lambda+1)}{\sqrt{\pi}~\Gamma(\lambda+1/2)}\int_{-1}^1
T_{2k}(x)C_n^{\lambda}(x)^2
(1-x^2)^{\lambda-1/2}dx\\
 &=\frac{(n+\lambda)~\Gamma(\lambda)~(2\lambda)_n}{n!~\Gamma(\lambda+1/2)}  \sum_{j=0}^n
\frac{(1/2-k)_k~(-j-\lambda)_k}{(1/2)_k~(-k-j-\lambda)_k}~
\frac{(-n)_j~(n+2\lambda)_j~(\lambda)_j}{(2\lambda)_j~(\lambda+1/2)_j~j!}
~\frac{\Gamma(j+\lambda+1/2)}{\Gamma(j+\lambda+1)}\\
&=(-1)^k~\frac{n+\lambda }{n!}\sum_{j=0}^n(-1)^j~
\left(\begin{array}{c}n\\j\end{array}\right)\frac{(-j-\lambda)_k}{(-k-j-\lambda)_k}~
(2\lambda+j)_n~\frac{1}{j+\lambda}\\
 &=\frac{n+\lambda}{n!}\sum_{j=0}^n(-1)^j~
\left(\begin{array}{c}n\\j\end{array}\right)
~\frac{(2\lambda+j)_n}{j+\lambda}
\frac{(-j-\lambda)_k}{(j+\lambda+1)_k},
\end{align*}
where we have used standard properties of the Gamma function and
Pochhammer's symbol. This proves (\ref{explicitm1}). Formula
(\ref{mkn}) follows easily from (\ref{explicitm1}) and the well
known relation $ \Gamma(x) \Gamma(1-x)=\pi/\sin(\pi x)$ applied
with $x=j+\lambda+1 >0$. In particular, it shows that
$|m_{2k+2,n}| \leq |m_{2k,n}|$; finally, (\ref{asymptoticsMkn}) is
a consequence of (\ref{mkn}) and Stirling formula.
\end{proof}

\medskip

In order to discuss the truncation error in
(\ref{series_for_entropu_Gegen}) we need to introduce the
following notation: fixed $n$, $\lambda>-1/2 $, $\lambda \notin
\N_0$, and $N \in \N$, $N>n+ \lambda$, let
\begin{equation*}
R_{n}^{\lambda}(N)=2 n
\sum_{k=N}^{\infty}\frac{c_{2k,n}^{\lambda}~m_{2k,n}^{\lambda}}{k}\,;
\end{equation*}
$|R_{n}^{\lambda}(N)|$ is the absolute error of approximation of
$E_n^\lambda $ if we truncate the series in
(\ref{series_for_entropu_Gegen}) after $k=N-1$. Proposition
\ref{prop:terminating} shows that it is convenient to take $N >
n+\lambda $. We consider the case $\lambda>0 $ (for negative
$\lambda$, see Remark at the end of this Section):
\begin{proposition} \label{prop:error}
Let $n \in \N$, $\lambda>0 $, $\lambda \notin \N$, and $N \in \N$,
$N > n+\lambda$. Then
\begin{equation}\label{errorProp}
 |R_{n}^{\lambda}(N)| \leq \mathcal{F}_{n}^{\lambda}(N):=
\frac{n (n+\lambda) }{N} \, \sum_{j=0}^n  \,
 \frac{1 }{(n-j)! \, j!} \,\frac{(2\lambda+j)_n } {j+\lambda}\,
\frac{|(-j-\lambda+1)_{N-1}|}{(j+\lambda+1)_{N-1}}\,.
\end{equation}
\end{proposition}

\begin{proof}
Let us denote
\begin{equation} \label{A}
A_{j, n}^\lambda= (-1)^j~ \binom{n}{j} \,
 \frac{n+\lambda} {n!} \frac{(2\lambda+j)_n}{j+\lambda} = (-1)^j  \,
 \frac{n+\lambda} {j+\lambda}\, \frac{(2\lambda+j)_n}{(n-j)! \, j!}.
\end{equation}
Then by (\ref{explicitm1}),
\begin{equation*}
R_{n}^{\lambda}(N)=2 n \, \sum_{k=N}^{\infty}m_{2k,n}^\lambda
\frac{c_{2k,n}^\lambda}{k} = 2n\, \sum_{j=0}^n A_{j, n}^\lambda
\sum_{k=N}^{\infty} \frac{(-j-\lambda)_k}{
(j+\lambda+1)_k}\frac{c_{2k,n}^\lambda}{k}\,,
\end{equation*}
so that
\begin{equation}
|R_{n}^{\lambda}(N)| \leq  2n \sum_{j=0}^n | A_{j, n}^\lambda
|\sum_{k=N}^{\infty}
\left|\frac{(-j-\lambda)_k}{(j+\lambda+1)_k}\right|~\left|\frac{c_{2k,n}^\lambda}{k}\right|
 \leq \frac{2 n}{N}  \sum_{j=0}^n |A_{j, n}^\lambda
|\sum_{k=N}^{\infty}
\left|\frac{(-j-\lambda)_k}{(j+\lambda+1)_k}\right|\,.
\label{error1}
\end{equation}
Taking into account that
$$
|A_{j, n}^\lambda|=(-1)^j A_{j, n}^\lambda, \quad \text{and} \quad
|(-j-\lambda)_k|=(-1)^{[j+\lambda-1] }~(-j-\lambda)_k, \text{ for
} k\geq n+[\lambda]+1,
$$
where $[\lambda ]$ stands for the integer value of $\lambda $, we
get from (\ref{error1}),
$$
|R_{n}^{\lambda}(N)|
\leq\frac{2n(-1)^{[\lambda+1]}}{N}\sum_{j=0}^n A_{j, n}^\lambda
\sum_{k=0}^{\infty}
\frac{(-j-\lambda)_{k+N}}{(j+\lambda+1)_{k+N}}\,.
$$
But
$$
\sum_{k=0}^{M}
\frac{(-j-\lambda)_{k+N}}{(j+\lambda+1)_{k+N}}=-\frac{1}{2}\,
\frac{(-j-\lambda+1)_{N-1} }{(j+\lambda+1)_{N-1} } +
\frac{(-j-\lambda+1)_{N+M}}{2\, (j+\lambda+1)_{N+M} } \,,
$$
so that
\begin{equation}\label{error2}
|R_{n}^{\lambda}(N)| \leq \frac{n(-1)^{[\lambda]}}{N}\sum_{j=0}^n
A_{j, n}^\lambda
\frac{(-j-\lambda+1)_{N-1}}{(j+\lambda+1)_{N-1}}=\mathcal{F}_{n}^{\lambda}(N)\,,
\end{equation}
and the statement follows.
\end{proof}

Given $\varepsilon >0$ we can use (\ref{errorProp}) in order to
find a (preferably, lowest) value $N_0 \in \N$ such that
$|R_{n}^{\lambda}(N)| \leq \varepsilon $. Obviously, the lower
bound for $N_0$ will be $N_0 = n+[\lambda]+1$. It is helpful to
get also an upper bound for such an $N_0$, that can be obtained
taking advantage of the geometric decay of $\F_n^\lambda (N)$. It
is based on the following
\begin{proposition}\label{prop:error2}  Let $n \in \N$, $\lambda>0 $, $\lambda \notin
\N$, and $N\in \N$, $N > n+\lambda$. Then for all $ h \in \N_0$,
\begin{equation}\label{error}
\F_{n}^{\lambda}(N+h)\leq \frac{F_n^\lambda (N)}{(N+\lambda+h)^{2
\lambda}}\,,
\end{equation}
where $\F_{n}^{\lambda}(N)$ is defined in (\ref{errorProp}), and
\begin{equation}\label{F}
F_n^\lambda (N) =|\lambda \sin (\pi \lambda ) |\,
\Gamma^2(\lambda)\, \left(\frac{ N-\lambda  }{  N+\lambda
}\right)^{N-\lambda -1}
 \frac{n\, e^{2 \lambda }}{\pi N}\,  \sum_{j=0}^n
 \frac{n+\lambda} {j+\lambda}\, \frac{(2\lambda+j)_n   }{(n-j)! \, j!}
 \, \frac{(\lambda)_j (\lambda+1)_j }{
(N-\lambda -j)_j \, (N+\lambda )_j}
\end{equation}
is a decreasing function in $N$, such that
$$
F_n^\lambda (N) =O\left(\frac{1}{N} \right)\,, \quad N \to \infty.
$$
\end{proposition}
\begin{proof}
We can use the identity
\begin{align*}
&\frac{(-j-\lambda+1)_{N+h-1}}{(j+\lambda+1)_{N+h-1}}=
\frac{(-j-\lambda+1)_j \, (\lambda+1)_j}{(N+h-\lambda-j)_j\,
(N+h+\lambda)_j}
\frac{(-\lambda+1)_{N+h-1}}{(\lambda+1)_{N+h-1}}\\
&=(-1)^j\, \frac{ (\lambda)_j (\lambda+1)_j}{(N+h-\lambda-j)_j\,
(N+h+\lambda)_j} \frac{(-\lambda+1)_{N-1}}{(\lambda+1)_{N-1}}
\frac{(N-\lambda)_h}{(N+\lambda)_h}\,;
\end{align*}
thus, with notation (\ref{A}) and by (\ref{error2}),
\begin{align}
|R_{n}^{\lambda}(N+h)| &\leq \frac{n(-1)^{[\lambda]}}{N}\,
\frac{(-\lambda+1)_{N-1}}{(\lambda+1)_{N-1}}
\frac{(N-\lambda)_h}{(N+\lambda)_h} \sum_{j=0}^n (-1)^j\, A_{j,
n}^\lambda
 \frac{ (\lambda)_j (\lambda+1)_j}{(N+h-\lambda-j)_j\,
(N+h+\lambda)_j} \nonumber  \\ & \leq
\frac{n(-1)^{[\lambda]}}{N}\,
\frac{(-\lambda+1)_{N-1}}{(\lambda+1)_{N-1}}
\frac{(N-\lambda)_h}{(N+\lambda)_h} \sum_{j=0}^n (-1)^j\, A_{j,
n}^\lambda
 \frac{ (\lambda)_j (\lambda+1)_j}{(N-\lambda-j)_j\,
(N+\lambda)_j} \nonumber \\ &=
\frac{(N-\lambda)_h}{(N+\lambda)_h}\,
 B_n^\lambda(N) \,, \label{cota1}
\end{align}
with
\begin{align}
 B_n^\lambda(N)=&\frac{n}{N}\,
\frac{|(-\lambda+1)_{N-1}|}{(\lambda+1)_{N-1}}
 \sum_{j=0}^n
 \frac{n+\lambda} {j+\lambda}\, \frac{(2\lambda+j)_n}{(n-j)! \, j!}
 \frac{ (\lambda)_j (\lambda+1)_j}{(N-\lambda-j)_j\,
(N+\lambda)_j} \nonumber \\
=&\frac{n}{N}\,
\frac{\Gamma(N-\lambda) \Gamma^2(\lambda) |\lambda \sin (\pi
\lambda )|}{\pi \Gamma(N+\lambda) }
 \sum_{j=0}^n
 \frac{n+\lambda} {j+\lambda}\, \frac{(2\lambda+j)_n}{(n-j)! \, j!}
 \frac{ (\lambda)_j (\lambda+1)_j}{(N-\lambda-j)_j\,
(N+\lambda)_j}\,, \label{balance}
\end{align}
where we have used again the identity $\Gamma(\lambda )
\Gamma(1-\lambda )\sin (\pi \lambda )=\pi$. Alternatively,
$B_n^\lambda(N)$ can be represented in terms of the following
truncating and balanced hypergeometric series, valid for $\lambda
\notin \Z$,
\begin{equation*}
B_n^\lambda(N)= \frac{1 }{N} \, \frac{|(-\lambda +1)_{N-1}
|}{(\lambda +1)_{N-1}}\, \frac{n+\lambda }{(n-1)!}\, \frac{(2
\lambda)_n }{\lambda }\,  _4F_3\left(\left.
 \begin{array}{cc}-n,~n+2\lambda,~\lambda,~\lambda\\
 2\lambda,-N+\lambda+1,N+\lambda\end{array}
 \right|1\right).
\end{equation*}
In order to simplify the expression of the error, we may use (cf.\
\cite[p.17]{Luke75}) that for $x\geq y \geq 1$,
\[
\frac{\Gamma(y) }{\Gamma(x)}\leq\frac{y^{y-1}~e^x}{x^{x-1}~e^y}\,.
\]
Hence,
$$
\frac{(N-\lambda)_h}{(N+\lambda)_h}=\frac{\Gamma(N+\lambda)}{\Gamma(N-\lambda)}\,
\frac{\Gamma(N-\lambda+h)}{\Gamma(N+\lambda+h)} \leq
\frac{\Gamma(N+\lambda)}{\Gamma(N-\lambda)}\, \left(\frac{
N-\lambda +h }{  N+\lambda +h }\right)^{N-\lambda +h-1} \frac{
e^{2 \lambda } }{ (N+\lambda +h )^{2 \lambda }  }\,.
$$
It is straightforward to verify that $ ( x/( x+t) )^{x-1} $ is
decreasing in $x$ for $x, t>0$, so that
$$
\frac{(N-\lambda)_h}{(N+\lambda)_h} \leq
\frac{\Gamma(N+\lambda)}{\Gamma(N-\lambda)}\, \left(\frac{
N-\lambda  }{  N+\lambda  }\right)^{N-\lambda -1} \frac{ e^{2
\lambda } }{ (N+\lambda +h )^{2 \lambda }  }\,.
$$
Gathering this inequality and (\ref{balance}) in (\ref{cota1}) we
obtain the statement of the Proposition.

\end{proof}

\begin{corollary}
Let $\varepsilon >0$, $n \in \N$, $\lambda>0 $, $\lambda \notin
\N$, and $N\in \N$, $N > n+\lambda$. Then if for $F_n^\lambda (N)$
defined in (\ref{F}),
\begin{equation} \label{trunc}
h \geq \max \left\{\left[ \left(\frac{F_n^\lambda
(N)}{\varepsilon}\right)^{\frac{1}{2\lambda}}-N-\lambda \right],
0\right\}\,,
\end{equation}
then
$$
\left|\F_{n}^{\lambda}(N+h)\right|\leq\varepsilon\,.
$$
\end{corollary}
Hence, if we want to find a suitable value of $N_0\in \N$ for
which $\left|R_{n}^{\lambda}(N_0)\right|\leq\varepsilon$, we can
use the following procedure:

\medskip

\begin{center}
\framebox{
\begin{minipage}{10cm}
\begin{center}
\textsc{Algorithm 3}
\end{center} \smallskip
\texttt{\begin{tabbing} \scriptsize{(1)} \phantom{1111} \= Take
$N=n+[\lambda ]+1$ and compute $ \F_{n}^{\lambda}(N)  $;
\\
\scriptsize{(2)} \> if $ \F_{n}^{\lambda}(N) \leq \varepsilon $,
put $N_0=N$ and quit;
\\
\scriptsize{(3)}  \> else compute $h$ equal to the r.h.s.\ of
(\ref{trunc});
\\
\scriptsize{(4)}  \> if $h> 0$, use bisection in order to find the
lowest $N_0 \in [N, N+h] \cap \N$ such that  \\  \> $
\F_{n}^{\lambda}(N_0) \leq \varepsilon $
\end{tabbing}}
\end{minipage}
}
\end{center}

\medskip

Obviously, we can use a more sophisticated zero-finding method in
the procedure above; however, usually bisection, which is simple
and easy to implement, is sufficient for our needs.

\medskip

\textsc{Remark: } A simple alternative for truncation of the
series in (\ref{series_for_entropu_Gegen}) can be computing
$E_n^\lambda$ for different integer values of $\lambda $ (say,
$[\lambda ]$ and $[\lambda ]+1$) and interpolating the value of
$E_n^\lambda$ for the non-integer  $\lambda $. Nevertheless, as
numerical experiments in Section \ref{sec:experiments} show, this
approach is not very satisfactory. For small values of $\lambda $
it yields large errors, and for large $\lambda $'s the loss in
speed computing the entropy at least twice can be compensated by
larger truncation values $N$.

\medskip

\textsc{Remark: } For small values of the parameter $\lambda $ the
error bounds above usually yield large truncation values $N$. This
occasionally might justify the use of explicit formulas
(\ref{explicitm1}) and (\ref{mkn}) for computation of $m_{2k,
n}^\lambda $ instead of Step 3. As $\lambda $ grows larger, the
explicit evaluation of Pochhammer's symbols rapidly becomes
substantially more time consuming and less accurate than matrix
multiplication.

Furthermore, for $-1/2 < \lambda <0$ the structure of the
coefficients $m_{2k, n}^{\lambda}$ yields extremely pessimistic
upper bounds for the error $|R_n^{\lambda}(N)|$, with a rate of
convergence even lower than $1/N$ established in Theorem
\ref{theoremMainNew}. Nevertheless, in this case the Gegenbauer
polynomials are uniformly bounded on $[-1,1]$, and truncation
error is  estimated better using formulas
(\ref{condSuf1New})--(\ref{remainderNew}).

\section{Numerical experiments for Gegenbauer polynomials} \label{sec:experiments}

\begin{figure}[hbt]
\centering \includegraphics[scale=0.8]{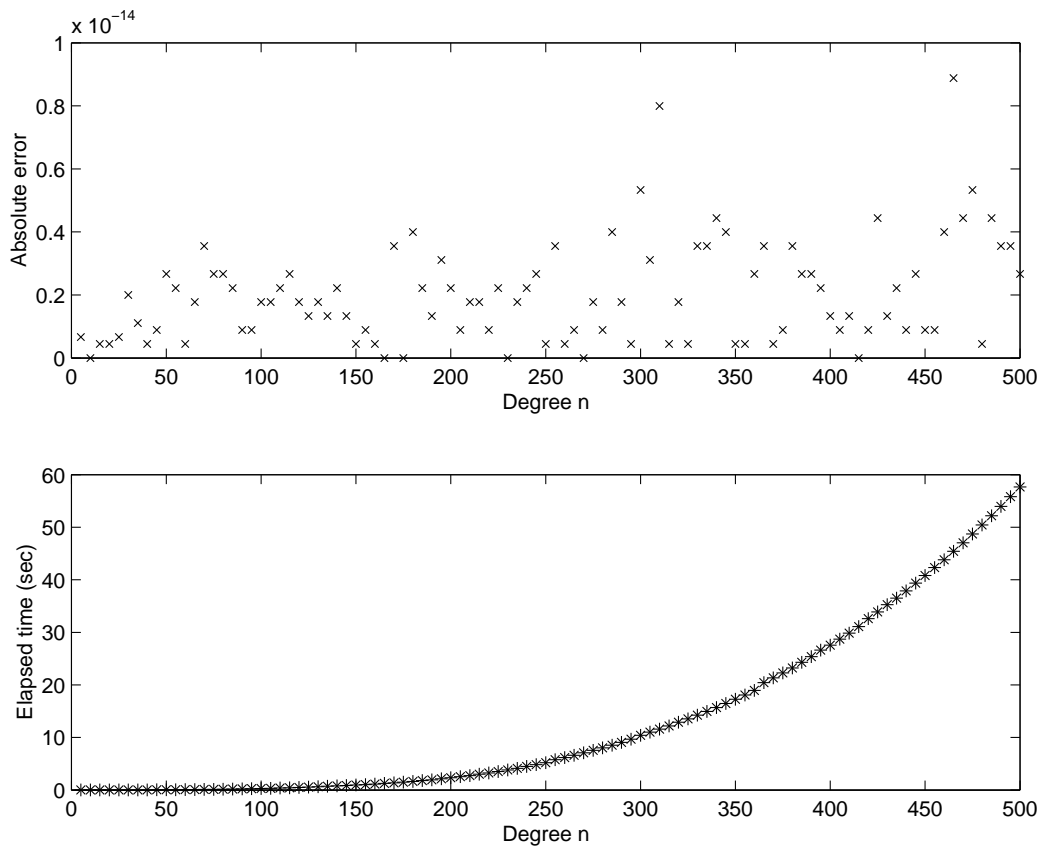}
\caption{Error and execution time of the algorithm for computation
of $E_n^2$.} \label{fig:Test_time_error_Gegen2}
\end{figure}

In this section we discuss briefly the performance of the
algorithm presented above, and compare it with some alternative
algorithms used for computing the entropy of Gegenbauer
polynomials.

We will compute $E_n ^\lambda$ for several values of $\lambda $.
In all cases Algorithm 2 was implemented in Matlab\texttrademark \
and executed on a computer with a single AMD Athlon\texttrademark
\ XP2000+ processor, 256 Mb RAM, and running Matlab 6 under
Windows. In particular, specific routines for sparse matrix
construction have been used through Matlab built in functions
\texttt{spdiags} and \texttt{speye}. For the experiments no
special performance optimization techniques have been used,
although we implemented vectorization when available.

One obvious test situation corresponds to $\lambda =0, 1, 2$, when
the explicit value of the entropy is known (cf.\ formulas
(\ref{ent1bis})--(\ref{ent2})).  Figure
\ref{fig:Test_time_error_Gegen2} shows that the error of the
algorithm (comparing with the exact value (\ref{ent2})) is
negligible. The execution time, computed as an average of 100 runs
of the algorithm, grows geometrically with the degree $n$ (due to
the proportional growth of the matrices involved), but still it is
below 1 minute for $n$ as large as 500.

\begin{figure}[hbt]
\centering \includegraphics[scale=0.85]{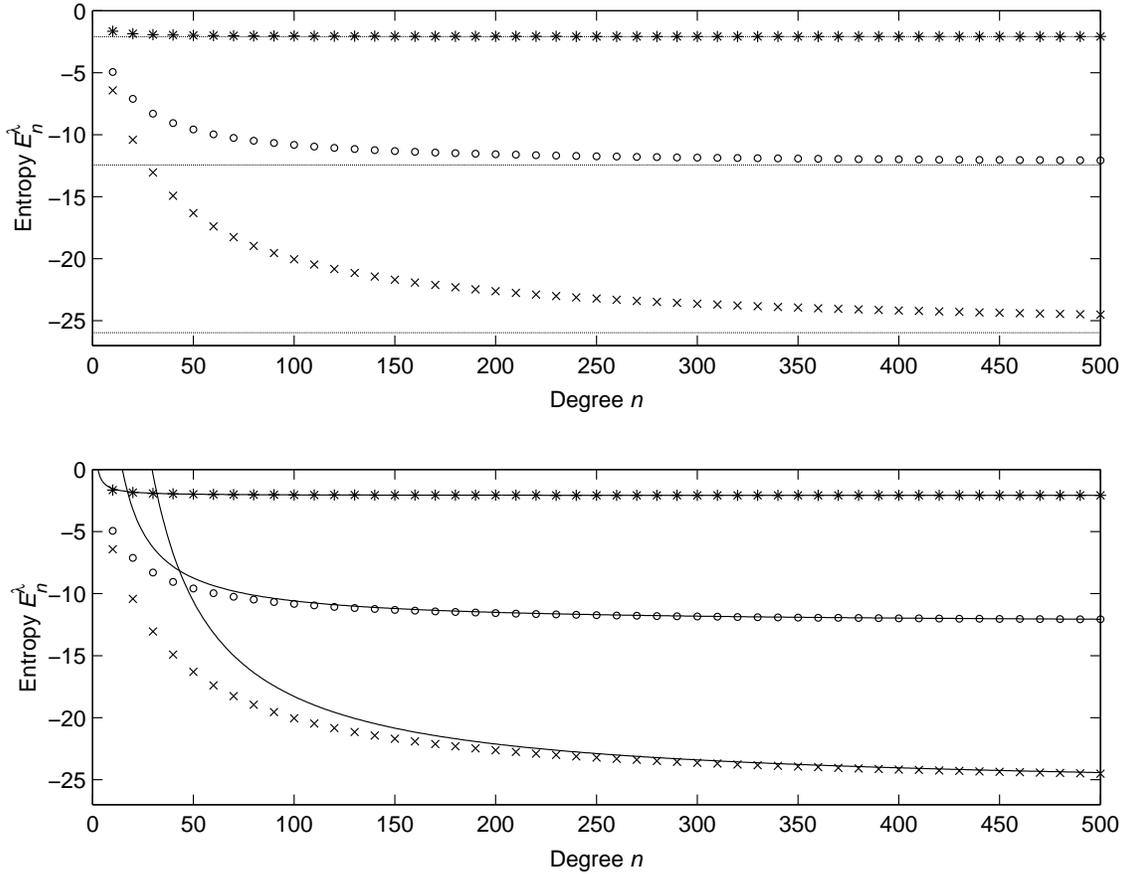}
\caption{Entropy $E_n^\lambda $ for $\lambda =2$ (marked with
'*'), $\lambda =10$ (marked with 'o') and $\lambda =20$ (marked
with 'x') compared with its limit $\mathcal{E}^\lambda_0$ (top)
and the asymptotic expression
$\mathcal{E}^\lambda_0+\mathcal{E}^\lambda_1/n$ (bottom).}
\label{fig:Test_asymp_Gegen}
\end{figure}

It is also illustrating to compare $E_n^\lambda $ with their
asymptotic expansion (\ref{asympt_Gegen_entropy}) truncated after
the first and the second terms (Fig.\ \ref{fig:Test_asymp_Gegen}).

In order to illustrate the performance of Algorithm 2 we compare
it with the following procedures to compute $E_n^\lambda $:
\begin{itemize}
    \item By formula
    (\ref{entropyGegenbDef}) using adaptative quadrature
    implemented in Mathematica\texttrademark \ 4.2;
    \item By formula
    (\ref{entropyGegenbDef}) using functions \texttt{quad} and \texttt{quadl} of Matlab 6,
    applied to explicit expressions of the polynomials $G_n^\lambda$ with coefficients
    computed using Mathematica 4.2 with exact arithmetics;
    \item By the algorithm described
    in \cite{MR1790053}. This approach is valid for integer values
    of $\lambda $ only.
\end{itemize}
In Table  \ref{tab1} we compare the errors and execution times of
the procedures above with those of Algorithm 2. The execution time
is taken as the average of 100 runs of the corresponding
algorithms.

\begin{table}[htbp]
\centering
\begin{tabular}{|c|c|c|c|c||c|c|c|c|} 
\hline \raisebox{-1.50ex}[0cm][0cm]{Method}  &
\multicolumn{4}{|c||}{Absolute error}& \multicolumn{4}{|c|}{Time (sec)} \\
\cline{2-9} & $ n=10\strut$ & $ n=25$& $ n=50$ & $ n=100$ & $
n=10$ &
$ n=25$& $ n=50$ & $n=100$ \\
 \hline (i)& $4.0 \times 10^{-7\strut}$\strut & $1.4
\times 10^{-4}$ & $5.8\times 10^{ -4}$ & $1.4  \times 10^{-3}$ &
$0.30\strut$& $0.55$ & $1.13$ &
$2.51$ \\
\hline (ii)& $4.3  \times 10^{-6\strut} $& $1.3 \times 10^{-5}$ &
$5.7  \times 10^{-5\strut} $ & $1.1  \times 10^{-2\strut} $ &
$0.08\strut$ & $0.28$& $0.89$ &
$2.61 $ \\
\hline (iii)& $7.7\times 10^{-8\strut}$ & $1.4\times 10^{-7}$ & $
8.6  \times 10^{-7\strut} $ & $5.2  \times 10^{-6\strut} $ &
$0.11\strut$& $0.31$& $0.74$&
$3.01$\\
\hline (iv)& $0 $& $0 $& $1.8 \times 10^{-26\strut}$ & $1.3 \times
10^{ -7}$& $0.006\strut$& $0.008$& $0.01$&
$0.02$  \\
\hline (v) & $4.4  \times 10^{-15\strut}$& $2.2   \times 10^{-16}$
& $6.0 \times 10^{-15}$ & $2.7 \times 10^{-15}$& $0.002$& $0.007$&
$0.03$&
$0.25$  \\
\hline
\end{tabular}

\medskip

\caption{Absolute error (left half) and time of computation (in
seconds) of $E_n^2$ for $n=10, 25, 50, 100$ by means of (i)
adaptative quadrature of Mathematica 4.1 with extended precision
(left half) or  using floating point arithmetics (right half),
(ii) \texttt{quad} and (iii) \texttt{quadl} functions of Matlab,
(iv) algorithm from \cite{MR1790053} implemented in Mathematica
4.1 with extended precision, and (v) Algorithm 2.} \label{tab1}
\end{table}

\begin{figure}[hbt]
\centering \includegraphics[scale=0.85]{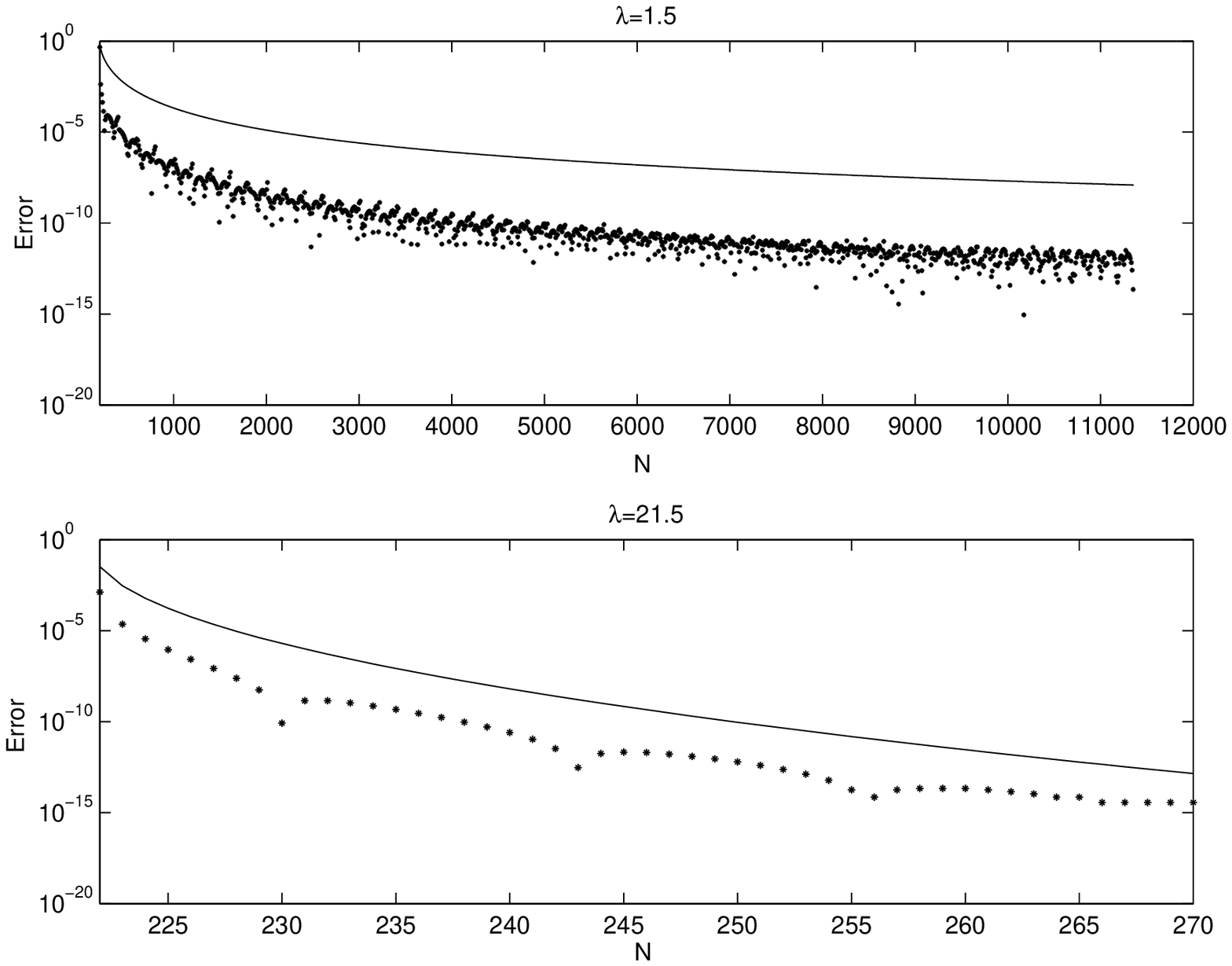} \caption{Error in
computing $E_{200}^\lambda $ by (\ref{series_for_entropu_Gegen})
truncating the series at $N$ (dots) compared with the error bound
$\mathcal{F}_{200}^{\lambda}(N)$ (solid line) for $\lambda =1.5$
(top) and $\lambda =21.5$ (bottom).} \label{fig:Test_F}
\end{figure}

As it was mentioned above, the bound in (\ref{errorProp}) usually
overestimates the error. For that purpose we truncate first the
series in (\ref{series_for_entropu_Gegen}) at $N$ such that
$\mathcal{F}_{n}^{\lambda}(N)$ is not greater than the machine
epsilon; the corresponding value of $E_n^\lambda $ is assumed as
the ``true'' value of the entropy. We compare it with the
approximation of $E_n^\lambda $ that we obtain if we truncate the
series in (\ref{series_for_entropu_Gegen}) at $N$ given by
Algorithm 3 (Fig.\ \ref{fig:Test_F}).

\begin{figure}[hbt]
\centering \includegraphics[scale=0.7]{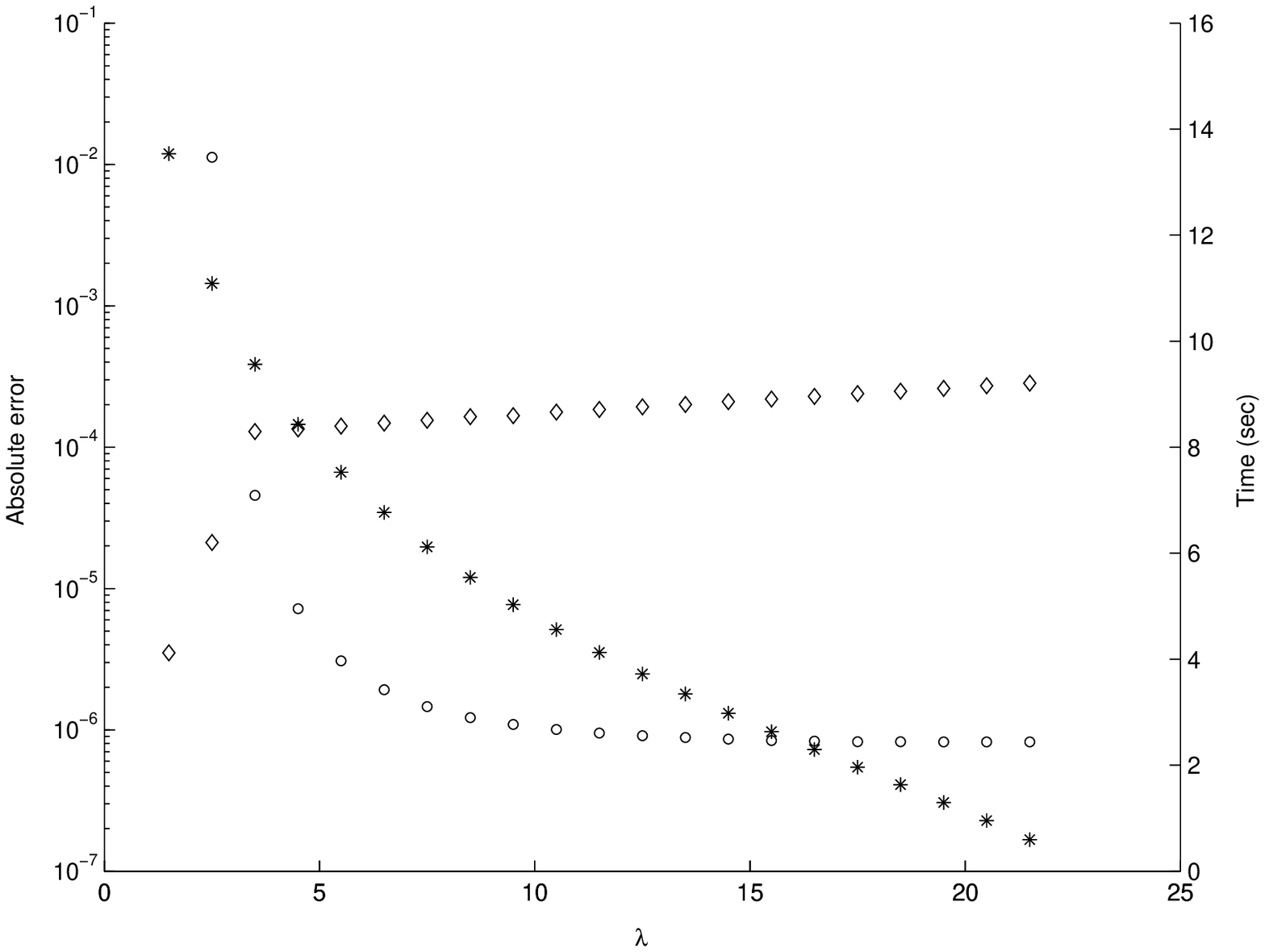}
\caption{Comparison of the interpolation strategy vs.\ truncation
of the series in (\ref{series_for_entropu_Gegen}).}
\label{fig:Test_interpol_Gegen.eps}
\end{figure}

Finally, it is tempting to avoid the question of truncation error
in (\ref{series_for_entropu_Gegen}) by applying Algorithm 2 to
$\lambda \in \N_0 $ and computing $E_n^\lambda $ for non-integer
values of the parameter by interpolation. As experiment, we
compute $E_{200}^{\lambda}$ for half-integer values of $\lambda$
using two strategies. In Figure \ref{fig:Test_interpol_Gegen.eps},
we observe the execution time (in seconds) for truncating the
series in (\ref{series_for_entropu_Gegen}) in the way that
$\F_{200}^{\lambda } \leq 10^{-6}$ (dots), and for interpolating
$E_{200}^{\lambda}$ by cubic splines using the values of the
entropy for $\lambda \pm 1/2$, $\lambda \pm 3/2$ (diamonds).
Asterisks represent the errors of interpolation. As we see,
interpolation usually yields large errors for small values of
$\lambda$, where it still could be competitive, since for large
$\lambda $'s we are penalized by the time invested in computing
the entropy at least twice.

\section{Computation of the entropy of spherical harmonics}
\label{sec:computationSH}

\begin{figure}[hbt]
\centering \includegraphics[scale=0.7]{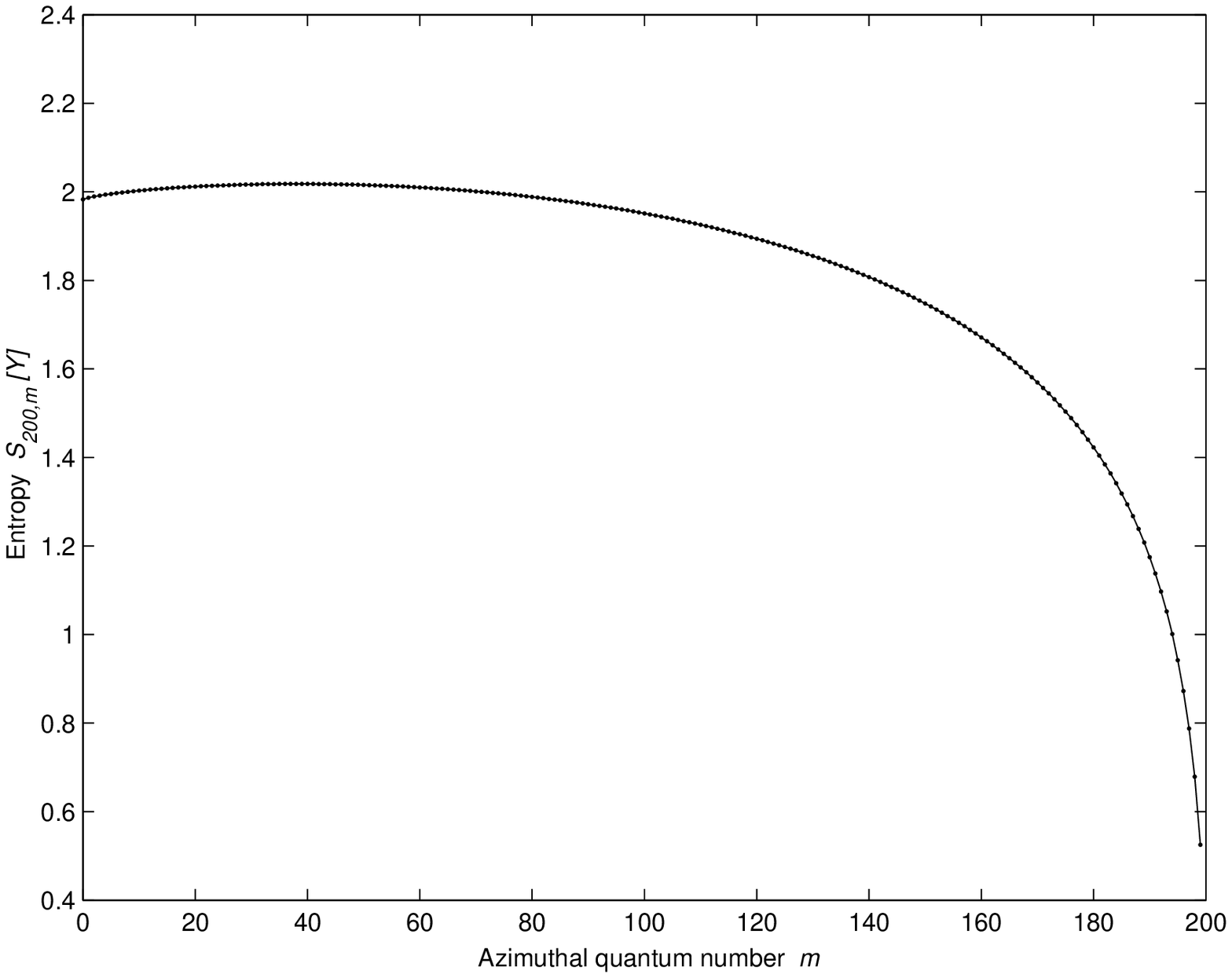}
\caption{Entropy $S_{200,m}[Y]$ for azimuthal quantum numbers $m
=0, 1,\dots, 199$.} \label{fig:Test_SphericalHarm}
\end{figure}

In this section, to show the usefulness of our computational
algorithm as well as the close connection of the entropy of
Gegenbauer polynomials analyzed in detail in the three previous
sections, we determine the spatial complexity of some
quantum-mechanical prototype and real systems with central
potentials (rigid rotator, harmonic oscillator, hydrogen atom,
Rydberg atoms, some diatomic molecules, etc.) by means of the
entropy of the spherical harmonics,
\[
S_{l,m}[Y] := - \int \left| Y_{l,m}(\Omega) \right| ^{2} \ln
\left| Y_{l,m}(\Omega) \right| ^{2} d \Omega,
\]
where $\Omega \equiv (\theta, \varphi)$, $d \Omega = \sin \theta
\, d \theta \, d \varphi$, with $0 \leq \theta \leq \pi$ and $0
\leq \varphi \leq 2\pi$, and $Y_{l,m}(\Omega)$ denotes the
spherical harmonics which depend on the orbital and azimuthal
quantum numbers, $l$ and $m$ respectively. It is well-known that
the principal quantum number $n \in \N_{0}$, together with $l \in
\N_{0}$ and $m \in \Z$, completely characterize a single-particle
system with a central potential. Moreover, for a given $n$ the
orbital quantum numbers $l \leq n-1$, and for a given $l$ the
azimuthal quantum number $-l \leq m \leq +l$.

The spatial or angular wavefunction of the system, which defines
its bulky shape, can be expressed in terms of Gegenbauer
polynomials as (cf.\ \cite{Yanez:99})
\[
Y_{lm}(\Omega) = N_{l,m}e^{im \varphi} (\sin \theta) ^{|m|} C_{l -
|m|}^{|m| + \frac{1}{2}}(\cos \theta),
\]
where $C_n^\lambda $ are Gegenbauer polynomials normalized as in
(\ref{standard}), and
\[
N_{l,m} = \left[ \frac{\left(l + \frac{1}{2} \right) \left( l -
|m| \right) ! \left[\Gamma \left( |m| + \frac{1}{2} \right)
\right] ^{2}}{2^{1-2 |m|} \pi ^{2} \left( l + |m| \right) !}
\right] ^{1/2}.
\]
Then, taking into account relation (\ref{normalizationGegen}), the
entropy $S_{l,m}[Y]$ is expressed in terms of the entropy of the
Gegenbauer orthonormal polynomials, defined in
(\ref{entropyGegenbDef}), as
\[
S_{l,m}[Y] = \ln \left(\frac{2\pi }{c_{|m|  + 1/2}}\right) +  E_{l
- |m|}^{|m| + 1/2}  - |m| \left[ 2\psi (l + |m| + 1) - 2\psi
\left( l + \frac{1}{2} \right) - 2 \ln 2 - \frac{1}{l+ 1/2}
\right]\,,
\]
where $c_\lambda $ is defined in (\ref{normalizationC}).

Thus, we can apply Algorithm 2 in order to compute the entropy of
the spherical harmonics $S_{l,m}[Y]$ for different values of the
quantum numbers $l$ and $m$. In Fig.\
\ref{fig:Test_SphericalHarm}, values of $S_{200,m}[Y]$ are
computed for integer values of $0\leq m \leq 200$. This figure
illustrates that for a given $l$ the entropy is higher around the
center of the manifold of azimuthal quantum numbers $m = -l,
-l+1,-l+2,\dots, l-1, l$, than at its extremes, indicating that
the spherical harmonics are much more localized for the largest
values of $|m|$. Moreover, the entropy is approximately constant
in the interval $-l/2\lesssim  m  \lesssim l/2$, and then it
monotonically decreases when $|m|$ grows up to its largest allowed
value $l$. The origin of this intriguing phenomenon is the
delicate interplay of the sinus factor and the Gegenbauer
polynomial involved in the spherical harmonics, which deserves
further numerical investigation.

\section*{Acknowledgements}

This work was partially supported by the Research Network on
Constructive Approximation (NeCCA), INTAS 03--51--6637 (V.S.B.,
J.S.D.\ and A.M.F.), a research grant from the Ministry of Science
and Technology (MCYT) of Spain, project code BFM2001-3878-C02
(J.S.D., A.M.F.\ and J.S.L.), and by Junta de Andaluc{\'\i}a, Grupo de
Investigaci{\'o}n FQM 0207 (J.S.D.) and FQM 0229 (A.M.F.\ and J.S.L.).
V.S.B.\ is indebted also to the Support Programm for Leading
Scientific Schools of Russia, project code NSh-1551.2003.1, and to
a research grant from the Russian Fund for Fundamental Research,
project code 02--01--00564.


\end{document}